\documentclass[letter, 12pt]{article}
        \usepackage{amsmath, amssymb, amsfonts, amsthm, latexsym}
        \usepackage[latin1]{inputenc}
        \usepackage{dsfont}
        \usepackage{amscd}
    \usepackage[all]{xy}

\newcommand{\ZZ}{\mathds Z}
\newcommand{\scrE}{\mathcal{E}}
\newcommand{\scrF}{\mathcal{F}}

\newcommand{\scrO}{\mathcal{O}}
\newcommand{\frcm}{\mathfrak{m}}

\newcommand{\lformula}[2]{\begin{equation} \label{#1} #2
\end{equation}}

\newcommand{\intlabel}[1]{\label{#1}}

\DeclareMathOperator{\Aut}{Aut}

\DeclareMathOperator{\cok}{cok}
\DeclareMathOperator{\Cov}{Cov}

\DeclareMathOperator{\Gal}{Gal}
\DeclareMathOperator{\Hom}{Hom}
\DeclareMathOperator{\id}{id}

\DeclareMathOperator{\Ind}{Ind}

\DeclareMathOperator{\rad}{rad}

\DeclareMathOperator{\Res}{Res}

\newtheorem{theorem}{Theorem} [section]
\newtheorem{lemma}[theorem]{Lemma}
\newtheorem{prop}[theorem]{Proposition}
\newtheorem{cor}[theorem]{Corollary}

\theoremstyle{definition}
\newtheorem{defn}[theorem]{Definition}

\theoremstyle{remark}

\begin{document}

\thispagestyle{empty} 
\centerline{\bf \LARGE Equivariant
Riemann-Roch theorems for} 
\centerline{\bf \LARGE curves over perfect fields}
\vskip5cm 
\centerline {\large 
\begin{tabular}{ll} Helena Fischbacher-Weitz & Bernhard K\"ock
\textit{(communicating author)}\\
Department of Pure Mathematics & School of Mathematics \\
University of Sheffield & University of Southampton \\
Hicks Building, Hounsfield Road & Highfield \\
Sheffield & Southampton \\
S3 7RH & SO17 1BJ \\
\texttt{h.fischbacher-weitz@shef.ac.uk} & \texttt{b.koeck@soton.ac.uk}
\\
 & Tel. 0044 23 8059 5125 \\
 & Fax 0044 23 8059 5147 
\end{tabular}}
\vskip5cm
\begin{abstract}
We prove an equivariant Riemann-Roch formula for divisors
on algebraic curves over perfect fields. 
By reduction to the known case of curves
over algebraically closed fields, we first show a preliminary formula
with coefficients in $\mathds{Q}$. 
We then prove and shed some further light on a divisibility result 
that yields a formula with integral coefficients. 
Moreover, we give variants of the main theorem for equivariant locally
free sheaves of higher rank.
\end{abstract}
\newpage
\setcounter{page}{1}

\section*{Introduction}
Let $X$ be a smooth, projective, geometrically irreducible curve
over a perfect field k and let $G$ be a finite subgroup of the
automorphism group $\Aut(X/k)$. For any locally free $G$-sheaf
$\scrE$ on $X$, we are interested in computing the equivariant Euler
characteristic
\[\chi(G,X,\scrE):=[H^0(X,\scrE)]-[H^1(X,\scrE)] \in K_0(G,k),\]
 considered as an
element of the Grothendieck group $K_0(G,k)$ of finitely generated
modules over the group ring $k[G]$. The main example of a locally
free $G$-sheaf we have in mind is the sheaf $\mathcal{L}(D)$
associated with a $G$-equivariant divisor $D = \sum_{P \in X}n_P P$
(that is $n_{\sigma(P)} = n_P$ for all $\sigma \in G$ and all $P \in
X$). If two $k[G]$-modules are in the same class in $K_{0}(G,k)$,
they are not necessarily isomorphic when the characteristic of $k$
divides the order of $G$. In order to be able to determine the
actual $k[G]$-isomorphism class of $H^0(X,\scrE)$ or $H^1(X,\scrE)$,
we are therefore also interested in deriving conditions for
$\chi(G,X,\scrE)$ to lie in the Grothendieck group $K_{0}(k[G])$ of
finitely generated \emph{projective} $k[G]$-modules and in computing
$\chi(G,X,\scrE)$ within $K_0(k[G])$.

The equivariant Riemann-Roch problem goes back to Chevalley and
Weil \cite{CW}, who described the $G$-structure of the space of global
holomorphic differentials on a compact Riemann surface. Ellingsrud
and L\o nsted \cite{EL} found a formula for the equivariant Euler
characteristic of an arbitrary $G$-sheaf on a curve over an
algebraically closed field of characteristic zero. Nakajima \cite{Na}
and Kani \cite{Ka} 
independently generalized this to curves over arbitrary
algebraically closed fields, under the assumption that the canonical
morphism $X \rightarrow X/G$ be tamely ramified. These results have
been revisited by Borne \cite{Bo}, who also found a formula 
that computes the difference between the equivariant Euler characteristics
of two $G$-sheaves in the case of a wildly ramified cover $X \rightarrow
X/G$. 
In the same setting, formulae for the equivariant
Euler characteristic of a single $G$-sheaf have been found by the
second author (\cite{BK1}, \cite{BK2}). Using these formulae, new
proofs for the reults of Ellingsrud-L\o nsted, Nakajima and Kani have
been given \cite{BK1}.

In this paper, we concentrate on the case
where the underlying field $k$ is perfect. Our main theorem, Theorem
3.4, is an equivariant Riemann-Roch formula in $K_0(k[G])$ when the
canonical morphism $X \rightarrow X/G$ is weakly ramified and
$\scrE=\mathcal{L}(D)$ for some equivariant divisor $D$. By
reduction to the known case of curves over algebraically closed
fields, we first show a preliminary formula with coefficients in
$\mathds{Q}$. The divisibility result needed to obtain a formula with
integral coefficients is then proved in two ways: Firstly, by applying
the preliminary formula to suitably chosen equivariant divisors; and
secondly, in two situations, by a local argument. The following
paragraphs describe the content of each section in more detail.

It is well-known that a finitely generated $k[G]$-module $M$ is
projective if and only if $M \otimes_k \bar{k}$ is a projective
$\bar{k}[G]$-module. In Section \ref{cartesian-d-section} we give 
a variant of this fact for classes in $K_0(G,k)$ rather than 
for $k[G]$-modules $M$ 
(Corollary \ref{Cartesian-d-individual-modules}). This
variant is much harder to prove and is an essential tool 
for the proof of our main result in Section \ref{strong-rr-formula-section}.

The first results in Section~\ref{strong-rr-formula-section} 
give both a sufficient condition and a necessary
condition under which the equivariant Euler characteristic 
$\chi(G,X,\scrE)$ lies in the image of the
Cartan homomorphism $c: K_0(G,k) \rightarrow K_0(k[G])$. More
precisely, when $\scrE = \mathcal{L}(D)$ for some equivariant
divisor $D = \sum_{P \in X} n_P P$, this holds if the canonical
projection $\pi: X \rightarrow X/G$ is weakly ramified and $n_P +1$
is divisible by the wild part $e_P^w$ of the ramification index
$e_P$ for all $P \in X$. When $\pi$ is weakly ramified we
furthermore derive from the corresponding result in \cite{BK2} the
existence of the so-called ramification module $N_{G,X}$, a certain
projective $k[G]$-module which embodies a global relation between
the (local) representations $\frcm_P/\frcm_P^2$ of the inertia group
$I_P$ for $P \in X$. If moreover $D$ is an equivariant divisor as
above, our main result, Theorem~\ref{my-theorem-2}, expresses
$\chi(G,X,\mathcal{L}(D))$ as an integral linear combination in
$K_0(k[G])$ of the classes of $N_{G,X}$, the regular representation
$k[G]$ and the projective $k[G]$-modules $\Ind_{G_P}^G(W_{P,d})$
(for $P\in X$ and $d\geq 0$) where the projective $k[G_P]$-module
$W_{P,d}$ is defined by the following isomorphism of
$k[G_P]$-modules:
\[\Ind_{I_P}^{G_P} (\Cov((\frcm_P/\frcm_P^2)^{\otimes (-d)})) \cong \bigoplus^{f_P} W_{P,d};\]
here Cov means taking the $k[I_P]$-projective cover and $f_P$
denotes the residual degree.
 
Finding an equivariant Riemann-Roch formula without denominators
amounts to showing that $W_{P,d}$ exists, i.e.\ that the left-hand side of
the above is ``divisible by $f_P$''. To do this, we use our prototype formula
\emph{with} denominators, formula (\ref{formula-with-fractions}), and
apply it to certain equivariant divisors $D$. 
If $\pi$ is tamely ramified, we
furthermore consider two situations where we can give a local proof
of the divisibility result, yielding a more concrete description of
$W_{P,d}$, see Proposition~\ref{structure-of-WP,d}. 

In Section 4, we give some variants of the main result that hold under
slightly different assumptions. In particular, these variants hold for
locally free $G$-sheaves that do not necessarily come from a divisor.

\section{Preliminaries}
\intlabel{prelim}
The purpose of this section is to fix some notations used throughout
this paper and to state some folklore results used later.

Throughout this section, let $X$ be a scheme of finite type over a
field $k$, and let $\bar k$ be an algebraic closure of $k$.
For any (closed) point $P \in X$, let $k(P):=\scrO_{X,P}/\frcm_{P}$
denote the residue field at $P$. 
Throughout this paper, let $\bar X$ denote the geometric fibre 
$X \times_k \bar k$, which is a scheme of finite type over $\bar k$, 
and let $p$
denote the canonical projection $\bar X \rightarrow X$.
Recall that $p$ is a closed, flat morphism which is in general 
not of finite type. We will see later that in dimension $1$, 
$p$ is ``unramified'' in
the sense that if $Q \in \bar X$ and $P=p(Q)$, then a local parameter
at $P$ is also a local parameter at $Q$. By Galois theory and
Hilbert's Nullstellensatz, we have for
every $P \in \bar X$:
\[ \# p^{-1} (P) = \# \Hom_k (k(P), \bar k) \leq
[k(P):k] < \infty,\] and equality holds if $k(P)/k$ is separable. 

Let now $G$ be a finite subgroup of $\Aut(X/k)$. Since the homomorphism
\[\Aut(X/k) \rightarrow \Aut(\bar X/
  \bar k), \sigma \mapsto \sigma \times \id\]
is injective, which is easy to check,
we may view $G$ as a subgroup of $\Aut(\bar X/\bar k)$.
Since the elements of $G$ act on the topological space of $X$ as
homeomorphisms, $G$ also acts on $|X|$, the set of closed points in
$X$. Analogously,
$G$ acts on the set $|\bar X|$ of closed points
in~$\bar X$.

\begin{defn} \intlabel{G-sheaf}
A \emph{locally free G-sheaf} (of rank $r$) on $X$
is a locally free $\scrO_{X}$-module $\scrE$ (of rank $r$) together
with an isomorphism of $\scrO_{X}$-modules $v_{\sigma}:\sigma^{*}\scrE
  \rightarrow \scrE$ for every $\sigma \in G$, such that for all
  $\sigma, \tau \in G$, the following diagram commutes:
\[\xymatrix{
\sigma^{*} \scrE \ar[r]^{v_{\sigma}} & \scrE\\
 \sigma^{*} (\tau^{*} \scrE)= (\tau \sigma)^{*} \scrE
 \ar[u]^{\sigma^{*} v_{\tau}} \ar[ur]_{v_{\tau \sigma}}}\]
\end{defn}
If $\scrE$ is a locally free $G$-sheaf
of finite rank, then the cohomology groups $H^{i}(X,\scrE)$
($i \in \mathds{N}_0$) are $k$-representations of $G$.
If moreover $X$ is proper over $k$, then the $H^{i}(X,\scrE)$
are finite-dimensional and vanish for $i>>0$
(see Theorem III.5.2 in \cite{Ha}).

We denote the Grothendieck group
of all finitely generated $k[G]$-modules
(i.e.\ finite-dimensional $k$-representations of $G$) by
$K_0(G,k)$, as opposed to the notation  $R_k(G)$
used by Serre in \cite{Se-g-f}.

\begin{defn} \intlabel{Euler-char-gen}
If $X$ is proper over $k$, and $\scrE$ is a
locally free $G$-sheaf of finite rank, then
\[\chi(G,X, \scrE):= \sum_{i} (-1)^{i} [H^{i}(X, \scrE)] \in K_0 (G,k)\]
is called the \emph{equivariant Euler characteristic} of $\scrE$ on $X$.
\end{defn}

For $P \in |X|$ or $P \in |\bar X|$, 
the \emph{decomposition group} $G_{P}$ and the \emph{inertia group}
$I_{P}$ are defined as follows:
\begin{align*}
G_{P}& :=\{\sigma \in G | \sigma(P)=P\} \text{;} \\
I_{P}& :=\{\sigma \in G_{P} | \bar \sigma = \id_{k(P)}\}=\ker(G_{P}
\rightarrow \Aut(k(P)/k)) \text{.} \end{align*}
Here $\bar \sigma$ denotes the endomorphism that $\sigma$
induces on $k(P)$.
Note that for all $Q \in |\bar X|$, we have $G_{Q}= I_{Q}$ and $G_{Q} = I_{P}$,
where $P:=p(Q) \in |X|$.

In the following lemma, we will assume for the first time that the
field $k$ is \emph{perfect}. 
\begin{lemma} \intlabel{tensoring-fibre}
Assume that $k$ is perfect.
Let $\scrF$ be a coherent sheaf on $X$, and let $\bar \scrF:=p^{*}
\scrF$. Let $P$ be a point in $X$, and let $\scrF(P):=\scrF_P
\otimes_{\scrO_{X,P}} k(P)$ be the fibre of
$\scrF$ at $P$.
Then the canonical homomorphism
\[\scrF(P) \otimes_{k} \bar k \mapsto \bigoplus_{Q \in p^{-1}(P)} \bar
\scrF (Q)\]
is an isomorphism.
In particular, the canonical homomorphism
\[k(P)\otimes_{k} \bar k \rightarrow \bigoplus_{Q \in p^{-1}(P)}k(Q)\]
is an isomorphism.
\end{lemma}
\begin{proof}
It follows from Galois theory that for any separable finite field
extension $k'/k$, the homomorphism
\[ k' \otimes_{k} \bar k \rightarrow \bigoplus_{\Hom_{k}(k',
  \bar k)} \bar k \]
defined by
\[y \otimes z \mapsto (\varphi(y) \cdot z)_{\varphi \in \Hom_{k}(k',
  \bar k)}\]
is an isomorphism. Since $k$ is perfect, by putting $k'=k(P)$ 
this implies the second part of the lemma, i.e. the special case 
where $\scrF=\scrO_X$.

Since the lemma is a local statement on $X$, we may assume that $X$ is
affine. The general case then follows from the special case together
with the definitions and basic properties of coherent sheaves and 
fibred products.
\end{proof}

\begin{prop} \intlabel{Omega(P)}
Assume that $k$ is perfect. Let $\Omega_{X/k}$ be the sheaf of
relative differentials of $X$ over
$k$. Then for every point $P \in |X|$, the canonical map
\[\frcm_{P}/\frcm_{P}^2 \rightarrow \Omega_{X/k}(P)\] is an isomorphism.
\end{prop}
\begin{proof}
Let $\Omega_{k(P)/k}$ denote the module
of relative differential forms of $k(P)$ over $k$. Using some
basic properties of differentials and of the cotangent space in an
affine setting, it
follows from Corollary 6.5 in \cite{Ku-K-d} that we have an exact sequence
\[0 \rightarrow \frcm_P/\frcm_P^2 \rightarrow
\Omega_{X/k}(P) \rightarrow
\Omega_{k(P)/k} \rightarrow 0. \]

By Corollary 5.3 in \cite{Ku-K-d}, $\Omega_{k(P)/k}$ is
trivial, so the map $\frcm_P/\frcm_P^2 \rightarrow
\Omega_{X/k}(P)$ is an isomorphism.
\end{proof}

Note that both Corollary 6.5 and Corollary 5.3 in \cite{Ku-K-d} require
$k(P)/k$ to be separable. 

Both Lemma \ref{tensoring-fibre} and Proposition \ref{Omega(P)} can be
turned into equivariant statements in the following sense. If we require
$\scrF$ to be a locally free $G$-sheaf, then for every point $P \in
|X|$, we obtain an action of the inertia group $I_{P}$ on the fibre
$\scrF(P)$ by $k(P)$-automorphisms. 

The action of $I_{P}$
on the fibre $\Omega_{X}(P)$ of the canonical sheaf corresponds to the
action on the cotangent space $\frcm_{P}/\frcm_{P}^2$ via the
isomorphism from Proposition \ref{Omega(P)}.

By letting $I_{P}$ act trivially on $\bar k$, we can extend the
action of $I_P$ on $\scrF(P)$ to an action on the tensor 
product $\scrF(P) \otimes_{k} \bar
k$. On the other hand, since $I_{Q}=I_{P}$ for any point $Q \in
p^{-1}(P)$, $I_{P}$ acts on
the fibre $\mathcal{G}(Q)$ of any locally free $G$-sheaf $\mathcal{G}$
on $\bar X$ for any point $Q \in p^{-1}(P)$. In particular, this holds
if $\mathcal{G}=p^{*}\scrF$ for a locally free $G$-sheaf $\scrF$ on~$X$.
With respect to these group actions, the isomorphism from Lemma
\ref{tensoring-fibre} is an isomorphism of $\bar k[I_{P}]$-modules.

We also have an action of the decomposition group $G_{P}$ on any fibre
$\scrF(P)$, but $G_{P}$ only acts on the fibre via $k$-automorphisms, whereas
$I_{P}$ acts via $k(P)$-automorphisms. $G_P$ does act
\emph{$k(P)$-semilinearly} on the fibre, that is,
for any $\sigma \in G_P, a \in
k(P)$ and $x,y \in \scrF(P)$ we have $\sigma.(ax + y)= (\bar \sigma.a)
(\sigma.x) + \sigma.y$, where $\bar \sigma$ denotes the automorphism of
$k(P)/k$ induced by $\sigma$.

Let now $X$ be a smooth, projective curve over a perfect field k.
Assume further  that $X$ is geometrically irreducible, i.e. that
the geometric fibre $\bar X=X \times_{k} \bar k$ is irreducible. 
Then the curve $X$ itself is irreducible.

The following lemma shows that although the canonical morphism $p:
\bar X \rightarrow X$ is usually not
of finite type, it can be thought of as an
``unramified'' morphism in the common sense, a fact that will be used
frequently throughout this paper.
\begin{lemma} \intlabel{p-unramified}
Let $Q \in |\bar X|$ be a closed point, and let $P:=p(Q)$.
Then every local parameter at $P$ is also a local parameter at $Q$.
\end{lemma}
\begin{proof}
Let $t_P$ be a local parameter at $P$. Then $t_P$ must be an element
of $\frcm_{P} \setminus \frcm_{P}^2$, so (the equivalence class of)
$t_{P}$ is a generator of the one-dimensional
vector space $\frcm_{P}/\frcm_{P}^2$ over $k(P)$.
Hence, $t_{P} \otimes 1$
is a generator of the rank-1 module
$\frcm_{P}/\frcm_{P}^2 \otimes_{k} \bar k$ over $k(P) \otimes_{k} \bar k$.

By Lemma \ref{tensoring-fibre} and Proposition \ref{Omega(P)}, we have a
canonical isomorphism
\[\frcm_{P}/\frcm_{P}^2 \otimes_{k} \bar k \rightarrow \bigoplus_{Q \in
  p^{-1}(P)} \frcm_{Q}/\frcm_{Q}^2\]
which we can view as an isomorphism of modules over $k(P) \otimes_{k}
\bar k \cong \bigoplus_{Q \in p^{-1}(P)} k(Q)$. 
Since this isomorphism must map $t_{P} \otimes 1$ to a
generator of the right-hand side over $\bigoplus_{Q \in
  p^{-1}(P)} k(Q)$, the image of $t_{P} \otimes 1$
in each component $\frcm_{Q}/\frcm_{Q}^2$ must be a generator of
$\frcm_{Q}/\frcm_{Q}^2$, i.e. the image of $t_{P}$ under each
induced homomorphism $p_{Q}: \scrO_{X,P} \rightarrow \scrO_{\bar
  X,Q}$ must be a local parameter at $Q$. 
\end{proof}
Let now $G$ be a finite subgroup of $\Aut(X/k)$.
It is a well-known result that the quotient scheme $Y:=X/G$ is also a
smooth projective curve, with function field $K(Y)=K(X)^{G}$. 
The canonical projection $X \rightarrow Y$ will be called $\pi$.
Let $P \in X$ be a closed point, $R:= \pi(P) \in Y$. 
Let $v_p$ be the unique normed valuation of the function field $K(X)$ 
associated to $P$, and let $v_R$ be the unique normed valuation of 
$K(Y)$ associated
to $R$. Then $v_P$ is equivalent to a valuation extending $v_R$. 
For $s \geq -1$,
we define the \emph{s-th ramification group $G_{P,s}$ at P} to be
the $s$-th ramification group
of the extension of local fields $K(X)_{v_{P}}/ K(Y)_{v_{R}}$. 
In particular, we have $G_{P,-1}=G_{P}$
  and $G_{P,0}=I_{P}$.

The canonical projection $\pi: X \rightarrow Y$ is called 
\emph{unramified (tamely ramified, weakly ramified)}  if $G_{P,s}$ is 
trivial for $s \geq 0 \; (s \geq 1, s \geq 2)$ and for all $P \in X$.
We denote the ramification index of $\pi$ at the place $P$
by $e_{P}$, its wild part by $e_{P}^{w}$ and its tame
part by $e_{P}^{t}$.
In other words, $e_{P}= v_{P}(t_{\pi(P)})= |G_{P,0}|$,
$e_{P}^{w}=|G_{P,1}|$ and $e_{P}^{t}=|G_{P,0}/G_{P,1}|$.

If $Q \in |\bar X|$ is a closed point, $P:=p(Q) \in |X|$, then for every
$s \geq 0$, we have $G_{Q,s}=G_{P,s}$ (by Proposition 5 in Chapter IV
in \cite{Se-c-l} and Lemma~\ref{p-unramified}).
In particular, we have
$e_{P}=e_{Q}$, $e_{P}^{w}=e_{Q}^{w}$ and $e_{P}^{t}=e_{Q}^{t}$.

\section{A Cartesian diagram of Grothendieck groups} \intlabel{cartesian-d-section}

A $k[G]$-module $M$ is projective if and only if $M \otimes_k \bar k$
is a projective $\bar k[G]$-module. In this section, we will now show 
variants of this well-known fact for classes in $K_0(G,k)$ rather than 
$k[G]$-modules.

Let $K_0(k[G])$ denote the Grothendieck group of finitely generated
projective $k[G]$-modules. This is a free group generated by the
isomorphism classes of indecomposable projective
$k[G]$-modules.
The Cartan homomorphisms $c: K_{0}(k[G])
\rightarrow K_{0}(G,k)$ and $\bar c: K_{0}(\bar k[G]) \rightarrow
K_{0}(G,\bar k)$ are injective (\cite{Se-g-f}, 16.1,
Corollary 1 of Theorem 35), so
$K_0(k[G])$ may be viewed as a subgroup of $K_0(G,k)$. The homomorphism
$$\beta: K_0(G,k) \rightarrow K_0(G, \bar k)$$ 
defined by tensoring with $\bar k$ over $k$ restricts to a homomorphism
$$\alpha: K_0 (k[G]) \rightarrow K_0(\bar k[G]).$$
By Proposition (16.22) in \cite{CR}, both
homomorphisms $\beta, \alpha$ are split injections.

\begin{prop} \intlabel{Cartesian-diagram}
The following diagram with injective arrows
is Cartesian, i.e. it commutes and viewing
the injections as inclusions, we have $K_{0}(\bar k[G]) \cap
K_{0}(G,k)= K_{0}(k[G])$.
\[\begin{CD}
K_{0}(k[G])  @>{\alpha}>> K_{0}(\bar k[G])\\
@V{c}VV         @VV{\bar c}V \\
K_{0}(G,k) @>>{\beta}> K_{0}(G, \bar k)
\end{CD}\]
\end{prop}
\begin{proof} The commutativity is obvious. Now consider the
  extended diagram (with exact rows)
\[\xymatrix{
0 \ar[r] & K_{0}(k[G]) \ar[r]^{\alpha} \ar[d]_{c} & K_{0}(\bar k[G])
\ar[r] \ar[d]_{\bar c} & M \ar[r] \ar[d]_{f} & 0 \\
0 \ar[r] & K_{0}(G,k) \ar[r]^{\beta} & K_{0}(G, \bar k)
\ar[r] & N \ar[r] & 0}\]
where $M= \cok \alpha$, $N= \cok \beta$, and $f$ is the homomorphism
  $M \rightarrow N$ induced by $\bar c$. By the Snake Lemma, there is an
  exact sequence of abelian groups
\[0 \rightarrow \ker c \rightarrow \ker \bar c \rightarrow \ker f
\rightarrow \cok c \text{,} \]
the first two modules being trivial since $c$ and $\bar c$ are injective.
Since $\alpha$ is a split injection, $M= \cok
\alpha$ is free over $\mathds{Z}$, and therefore $\ker f$ must also be free
over $\mathds{Z}$. On the other hand, by Theorem (21.22) in \cite{CR},
we have $|G| \cdot \cok c =0$, so $\cok c$ is a torsion
module. Using the exactness of the sequence above, this implies
$\ker f=0$. Now an easy diagram chase completes the proof.
\end{proof}

Proposition \ref{Cartesian-diagram} 
says that given a class $\mathcal{C}$ in $K_0(G,k)$,
$\mathcal{C}$ lies in the image of $c$ if and only if
$\beta(\mathcal{C})$ lies in the image of $\bar c$. The following
corollary appears to be only slightly different from this, 
yet some additional tools will
be required for its proof.

\begin{cor} \intlabel{Cartesian-d-individual-modules}
Let $\mathcal{C}$ be a class in $K_{0}(G,k)$. Then $\mathcal{C}$ is
the class of a projective $k[G]$-module if and only if $\beta
(\mathcal{C})$ is the class of a projective $\bar k[G]$-module.
\end{cor}
Before proving Corollary \ref{Cartesian-d-individual-modules}, we will
need a few preliminary results on $k[G]$-modules.
Recall that a $k[G]$-module is called \emph{simple} if it is nonzero
and has no proper $k[G]$-submodules, and \emph{indecomposable} if it
is nonzero and is not a direct sum of proper $k[G]$-submodules.
\begin{prop} \intlabel{tensor-and-decompose}
\begin{enumerate}
\item[(a)] For every simple
  $k[G]$-module $M$, the $\bar k[G]$-module $M\otimes_{k} \bar k$ is
  semisimple.
\item[(b)]
Let $\{P_{1}, \ldots , P_{s}\}$ be a set of representatives of the
isomorphism classes of indecomposable projective
  $k[G]$-modules, and let
\[ P_{i} \otimes_{k} \bar k = \bigoplus_{j=1}^{r_{i}} \bar Q_{ij},
\;\; \bar Q_{ij}
  \;\text{indecomposable projective} \; \bar k[G] \text{-modules.} \]
Then every indecomposable $\bar k[G]$-module is isomorphic to some
$\bar Q_{ij}$. Further $\bar Q_{ij}~\cong~\bar Q_{i'j'}$ implies that $i=i'$,
i.e. there is no overlap between the sets of indecomposable $\bar k[G]$-modules
which come from different indecomposable $k[G]$-modules.
\end{enumerate}
\end{prop}

\begin{proof}
This proposition is a variation of Theorem 7.9 in \cite{CR}. In
\cite{CR}, the algebraic closure $\bar k$ is replaced by a finite
algebraic extension $E$ of $k$, and part (b) is stated for
\emph{simple} modules rather than for indecomposable projective
modules. Using only elementary algebraic methods, it can be shown
that there is a finite algebraic extension $E/k$ such that every
simple $\bar k[G]$-module can be realized as a simple $E[G]$-module,
i.e. every simple $\bar k[G]$-module $M$ can be written as $M=N
\otimes_{E} \bar k$ for some simple $E[G]$-module $N$. This suffices
to derive part~(a) from the result in \cite{CR}. Furthermore, it is
well-known that mapping every projective $k[G]$-module $P$ to the
$k[G]$-module $P/\rad P$ gives a 1-1 correspondence between the
isomorphism classes of indecomposable projective $k[G]$-modules and
the isomorphism classes of simple $k[G]$-modules, whose inverse is
given by taking $k[G]$-projective covers. We can thus deduce our
proposition from the result in \cite{CR}, using that projective
covers are additive (by Corollary 6.25 (ii) in \cite{CR}) and
commute with tensor products (by Corollary 6.25 (i) in \cite{CR}).
\end{proof}

\begin{proof}[Proof of Corollary \ref{Cartesian-d-individual-modules}]
The ``only if'' direction is obvious. For the ``if'' direction, we
note first of all that if $\mathcal{C}$ is a class in $K_0(G,k)$ and
$\beta(\mathcal{C})$ is the class of a projective $\bar k[G]$-module,
then Proposition \ref{Cartesian-diagram} yields that $\mathcal{C}$
can be viewed
as a class in $K_0(k[G])$. Hence it suffices to show the ``if''
direction for classes $\mathcal{C} \in K_0 (k[G])$, replacing the
homomorphism $\beta$ by its restriction $\alpha$.

Let $\{P_1, \ldots , P_{s}\}$  be a set of representatives of the
isomorphism classes of indecomposable $k[G]$-modules. 
Every $\mathcal{C} \in K_0(k[G])$ can now be written as a $\mathds{Z}$-linear
combination of
the classes $[P_i]$, and all coefficients of this linear combination
are nonnegative if and only if $\mathcal{C}$ is the class of a 
projective module.
Using Proposition \ref{tensor-and-decompose}, one now
easily shows that if $\alpha(\mathcal{C})$ is the class of a
projective module in $K_0(\bar k[G])$, then $\mathcal{C}$ is the class of a
projective module in $K_0(k[G])$, which proves the assertion.
\end{proof}

\section{The equivariant Euler characteristic in terms of
  projective $k[G]$-modules} \intlabel{strong-rr-formula-section}
By a theorem of Nakajima, the equivariant Euler
characteristic of any locally free $G$-sheaf on $X$ lies in the
image of the Cartan homomorphism $c~:~K_0(k[G])~\rightarrow~
K_0(G,k)$, provided that the canonical projection $\pi: X
\rightarrow Y=X/G$ is \emph{tamely ramified}. In this section, we will
also consider the more general
case where $\pi$ is \emph{weakly ramified}. We give both a necessary
condition and a sufficient condition for the 
equivariant Euler characteristic to
lie in the image of $c$, provided that the $G$-sheaf in question has
rank 1 (comes from a divisor). Under this condition, we state an
equivariant Riemann-Roch formula in the Grothendieck group of
projective $k[G]$-modules.

We make the same assumptions and use the same notations as in
section \ref{prelim}. In particular $p$ denotes the
projection $\bar X= X \times_k \bar k \rightarrow X$. Additionally,
let $\bar \pi$ denote the canonical projection $\bar X \rightarrow
\bar Y:=\bar X/G= Y \otimes_k \bar k$, and let $\tilde p$ denote the
projection $\bar Y \rightarrow Y$. We have the following commutative
diagram:
\[\begin{CD}
\bar X @>{p}>> X\\
@V{\bar \pi}VV         @VV{\pi}V \\
\bar Y @>{\tilde p}>> Y
\end{CD}\]

\begin{theorem} \intlabel{Euler-char-projective-tam}
If $\pi$ is tamely ramified
and $\scrE$ is a locally free $G$-sheaf on
$X$, then the equivariant Euler characteristic $\chi(G,X, \scrE)$ lies in the
image of the Cartan homomorphism $c: K_{0}(k[G]) \rightarrow
K_{0}(G,k)$.
\end{theorem}
\begin{proof} Follows directly from Theorem 1 in \cite{Na}.
\end{proof}
\begin{theorem} \intlabel{Euler-char-projective}
Let $D=\sum_{P \in |X|} n_{P} P$ be a $G$-equivariant divisor on $X$.
\begin{enumerate}
\item[(a)] If $\pi$ is weakly ramified and $n_{P} \equiv -1 \mod e_{P}^{w}$
  for all $P \in X$, then the equivariant Euler characteristic $\chi(G,X, \mathcal{L}(D))$ lies in the
image of the Cartan homomorphism $c: K_{0}(k[G]) \rightarrow
K_{0}(G,k)$. If moreover one of the cohomology groups $H^{i}(X,
  \mathcal{L}(D))$, $i=0,1$, vanishes, then the other one is a
  projective $k[G]$-module.
\item[(b)] Let $\deg D > 2 g_{X} - 2$. If the $k[G]$-module $H^{0}(X,
  \mathcal{L}(D))$ is projective, then $\pi$ is weakly ramified and
  $n_{P} \equiv -1 \mod e_{P}^{w}$ for all $P \in |X|$.
\end{enumerate}
\end{theorem}
\begin{proof}
If $k$ is algebraically closed, the theorem coincides with Theorem
2.1 in \cite{BK2}. 

In the general case, if $\pi$ is weakly ramified and $D$ 
satisfies the congruence condition
``$n_P \equiv -1 \mod e_P^w $ for all $P$'', then $\bar
\pi: \bar X \rightarrow \bar Y$ is weakly ramified, and by 
Lemma~\ref{p-unramified}, the divisor $p^*D$ on $\bar X$ also satisfies the 
congruence condition. 
By the special case, $\chi(G,X,\mathcal{L}(p^*D)$ then lies in the image of 
$\bar c$. Hence by
Proposition \ref{Cartesian-diagram}, 
$\chi(G,X,\mathcal{L}(D)$ lies in the image of
$c$. Here we have used that
$H^i(X, \mathcal{L}(D)) \otimes_k \bar k = H^i(\bar X, \mathcal{L}(p^*
D))$ for every $i$ (cf. Proposition III.9.3 in \cite{Ha}). This
also implies the rest of part (a).

For part (b), let $\deg D > 2 g_X -2$. and let $H^0(X,
\mathcal{L}(D))$ be projective. Then $\deg p^* D > 2 g_{\bar X} -2$
and $H^0(\bar X, \mathcal{L}(D))$ is
projective. Thus $\bar \pi: \bar X \rightarrow \bar Y$ is weakly
ramified and the congruence condition holds. But then $\pi$ is
weakly ramified also, and the congruence condition holds for $D$,
again by Lemma~\ref{p-unramified}.
\end{proof}

The following theorem generalizes Theorem 4.3 in \cite{BK2} and will
be used in the formulation of the (main) Theorem~\ref{my-theorem-2}. 
We refer the
reader to page~1101 of the paper \cite{BK2} for an account of the
nature, significance and history of the ``ramification module''
$N_{G,X}$ and for simplifications of formulae (\ref{NG,X-first-eqn})  
and (\ref{NG,X-second-eqn})  when
$\pi$ is tamely ramified.
\begin{theorem} \intlabel{NG,X}
Let $\pi$ be weakly ramified. Then there is a projective
$k[G]$-module $N_{G,X}$ such that
\lformula{NG,X-first-eqn}{\bigoplus^{n} N_{G,X} \cong \bigoplus_{P
\in X} \bigoplus_{d=1}^{e_{P}^{t}-1} \bigoplus^{e_{P}^{w} \cdot  d}
\Ind_{I_{P}}^{G}(\Cov((\frcm_P/\frcm_P^2)^{\otimes d})) \text{,}}
where $\Cov$ denotes the $k[I_{P}]$-projective cover. The class of
$N_{G,X}$ in $K_{0}(G,k)$ is given by
\lformula{NG,X-second-eqn}{[N_{G,X}]= (1-g_{Y}) [k[G]] - \chi(G,X,
\mathcal{L}(E))} where $E$ denotes the $G$-equivariant divisor
$E:=\sum_{P \in X}
  (e_{P}^{w} -1) \cdot P$.
\end{theorem}
\begin{proof}
Theorem 4.3 in \cite{BK2} yields that there is a projective $\bar
k[G]$-module $N_{G,\bar X}$ such that 
\[\bigoplus^{n} N_{G,\bar X} \cong \bigoplus_{Q \in \bar X}
\bigoplus_{d=1}^{e_{Q}^{t}-1} \bigoplus^{e_{Q}^{w} \cdot d}
\Ind_{G_{Q}}^{G}(\Cov((\frcm_Q/\frcm_Q^2)^{\otimes d})) \text{,}\]
and that the class of $N_{G,\bar X}$ is given by
\[ [N_{G,\bar X}] =(1- g_{\bar Y}) [\bar
k[G]] - \chi(G,X, \mathcal{L}(\bar E))\] where 
$\bar E:=\sum_{Q \in \bar X} (e_{Q}^{w}-1)
\cdot Q = p^*E$.
Thus $[N_{G,\bar X}] = \beta (\mathcal{C})$ where
\[\mathcal{C}:= (1-g_{Y}) [k[G]] - \chi(G,X,\mathcal{L}(E)) \in K_{0}(G,k).\]
By Corollary \ref{Cartesian-d-individual-modules}, $\mathcal{C}$ is
the class  of
some projective $k[G]$-module, say $N_{G,X}$.
Using Lemma \ref{tensoring-fibre} and the injectivity of $\beta$, one 
easily shows that $N_{G,X}$ satisfies Formula (\ref{NG,X-first-eqn}).
\end{proof}

For every point $P \in X$,
let $f_{P}$ denote the residual degree $[k(P):k(\pi(P))]$.

\newpage
\begin{theorem}[Equivariant Riemann-Roch formula]\intlabel{my-theorem-2}
Let $\pi$ be weakly ramified.
\begin{enumerate}
\item[(a)] Let $P \in |X|$ be a closed point.
For every $d \in \{0, \ldots, e_P^t -1\}$, there is a unique
projective $k[G_P]$-module $W_{P,d}$ such that
\[\Ind_{I_P}^{G_P} (\Cov((\frcm_P/\frcm_P^2)^{\otimes (- d)})) \cong \bigoplus^{f_P} W_{P,d}\]
as $k[G_P]$-modules.

\item[(b)] Let $D= \sum_{P \in X} n_{P} \cdot P$ be a divisor on
$X$ with $n_{P}
\equiv -1 \mod e_{P}^{w}$ for all $P \in X$. For any $P
\in X$, we write
\[ n_{P}= (e_{P}^{w} -1)+ (l_{P} + m_{P} e_{P}^{t}) e_{P}^{w} \]
with $l_{P} \in \{0, \ldots , e_{P}^{t}-1\}$ and $m_{P} \in
\mathds{Z}$. Furthermore, for any $R \in Y$, fix a point $\tilde R
\in \pi^{-1}(R)$. Then we have in $K_{0}(k[G])_{\mathds{Q}}$:
\begin{multline} \label{Euler-char-divisor} 
\chi(G,X, \mathcal{L}(D)) \\
 = -[N_{G,X}] + \sum_{R \in Y} \sum_{d=1}^{l_{\tilde
R}}[\Ind_{G_P}^G(W_{P,d})] + \left(1 - g_{Y} + \sum_{R \in Y}
[k(R):k] m_{\tilde R} \right) [k[G]]. \end{multline}
\end{enumerate}
\end{theorem}
\begin{proof}
We first show that under the preconditions of (b), the following holds
in the Grothendieck group with rational coefficients $K_0(k[G])_{\mathds{Q}}$:
\begin{multline} \label{formula-with-fractions}
\chi(G,X, \mathcal{L}(D)) 
= - [N_{G,X}] + \sum_{R \in Y} \frac{1}
{f_{\tilde R}} \sum_{d=1}^{l_{\tilde R}}[\Ind_{I_{\tilde R}}^{G}
 (\Cov((\frcm_{\tilde R}/\frcm_{\tilde R}^2)^{\otimes (- d)}))] \\ 
+ \Bigm(1 - g_{Y} + \sum_{R \in Y} [k(R):k]
m_{\tilde R} \Bigm) [k[G]] \end{multline} 
With suitably chosen
divisors $D$, Formula (\ref{formula-with-fractions}) 
will then be used to show part (a). Formula
(\ref{formula-with-fractions}) and part (a) 
obviously imply part (b).

For curves over algebraically closed fields, we have $f_P=1$
for all $P$, so Formula (\ref{formula-with-fractions}) 
coincides with Theorem 4.5 in \cite{BK2}. 

The injective homomorphism $\beta:~K_0(G,k)~\rightarrow~K_0(G, \bar k)$ maps
$\chi(G,X,\scrE)$ to $\chi(G, \bar X, p^{*} \scrE)$, and by Theorem
\ref{Euler-char-projective}, both of these Euler characteristics lie
in the image of the respective Cartan homomorphisms. Hence  it suffices
to show that $\beta$ maps every summand of the right-hand side of
formula (\ref{Euler-char-divisor}) (applied to $X, D$) to the
corresponding summand of the right-hand side applied to $\bar X,
p^{*}D$.

From the proof of Theorem \ref{NG,X}, we see that
$\beta([N_{G,X}])=[N_{G, \bar X}]$.

By Lemma \ref{p-unramified}, we have $l_{Q}=l_{P}$ and
$m_{Q}=m_{P}$ whenever $Q \in p^{-1}(P)$. Furthermore, the number of
preimages of a point $R \in Y$ under $\pi: X \rightarrow Y$ 
is $\frac{n}{e_{\tilde R} f_{\tilde
R}}$. For any $S \in |\bar Y|$,
fix a point $\tilde S \in \bar \pi^{-1}(S)$.
Using Lemma \ref{tensoring-fibre}, we see that
\begin{align*} \beta & \left(\sum_{R \in Y} \frac{1}{f_{\tilde R}}
\sum_{d=1}^{l_{\tilde R}}[\Ind_{I_{\tilde R}}^{G}
 (\Cov((\frcm_{\tilde R}/ \frcm_{\tilde R}^2)^{\otimes (-d)}))] \right) \\
&  = \sum_{Q \in \bar X} \frac{e_{Q}}{n}
 \sum_{d=1}^{l_{Q}} [\Ind_{G_{Q}}^{G}
 (\Cov((\frcm_Q/\frcm_Q^2)^{\otimes (-d)}))]  \\
& = \sum_{S \in \bar Y} \sum_{d=1}^{l_{\tilde S}} [\Ind_{G_{\tilde
 S}}^{G} (\Cov((\frcm_{\tilde S}/ \frcm_{\tilde S}^2)^{\otimes (-
d)}))]\\
\end{align*}
Moreover, we have
\[\beta  \left( \bigm( 1 - g_{Y} + \sum_{R \in Y} [k(R):k] m_{\tilde
  R} \bigm)\; [k[G]] \right) = \left(1 - g_{\bar Y} + \sum_{S \in \bar Y} m_{\tilde S}\right)
[\bar k[G]]\text{,}\]
which completes the proof of Formula (\ref{formula-with-fractions}).

We now prove part~(a). Let $P \in X$ be a closed point. For $d=0$,
the statement is obvious because $(\frcm_P/\frcm_P^2)^0$ is 
the trivial one-dimensional $k(P)$-representation of $I_P$, so it
decomposes into $f_P$ copies of the trivial one-dimensional 
$k(R)$-representation of $I_P$, where $R:=\pi(P)$.
Hence we only need to do the inductive step from $d$ to $d+1$, for
$d \in \{0, \ldots, e_P^t -2 \}$.

If $\pi$ is unramified at $P$, then $e_P^t=1$, so there is no $d \in
\{0, \ldots, e_P^t -2 \}$. Hence we may assume that $\pi$ is
ramified at $P$. Set $H:=G_P$, the decomposition group at $P$, and let
$\pi'$ denote the projection $X \rightarrow X/H =:Y'$. For every
closed point $Q \in |X|$ and for every $s \geq -1$, 
let $H_{Q,s}$ be the $s$-th ramification group at $Q$ with respect to
that cover, as introduced in Section \ref{prelim}. Then we have 
$H_{Q, s}=G_P \cap G_{Q,s}$ for every $s \geq -1$ and
every $Q \in |X|$. In particular, if $\pi$ is weakly ramified, then
so is $\pi'$. For $Q=P$, we get $H_{P,s}=G_{P,s}$ for all $s \geq
-1$; in particular, the ramification indices and residual degrees of
$\pi$ and $\pi'$ at $P$ are equal.

Let now $D:= \sum_{Q \in |X|} n_Q \cdot Q$ be the
  $H$-equivariant divisor with coefficients
\[ n_{Q}= \left\{ \begin{array}{rl}
(d+2) e_Q^w -1 & \text{if} \quad  Q=P \\
e_Q^w -1 &  \text{otherwise} \end{array} \right.\]
Then formula (\ref{formula-with-fractions}) applied to $H, X, D$ gives
\begin{multline}
\chi(H,X, \mathcal{L}(D)) = - [N_{H,X}] + \frac{1}{f_P}
\sum_{n=1}^{d}[\Ind_{I_P}^{H}
 (\Cov((\frcm_P/\frcm_P^2)^{\otimes (- n)}))] \\
+ \frac{1}{f_P} [\Ind_{I_P}^{H}
 (\Cov((\frcm_P/\frcm_P^2)^{\otimes (- (d+1))}))] +
 (1 - g_{Y'})  [k[H]] \end{multline}
in $K_0(k[H])_\mathbb{Q}$.  By the induction hypothesis, the sum
from $n=1$ to $d$ in this formula is divisible by $f_P$ in
$K_0(k[H])$; hence the remaining fractional term
$\frac{1}{f_P}[\Ind_{I_P}^{H}
 (\Cov((\frcm_P/\frcm_P^2)^{\otimes (- (d+1))}))]$ must lie in
$K_0(k[H])$. In other words, when writing
$\Ind_{I_P}^H(\Cov((\frcm_P/\frcm_P^2)^{\otimes (-(d+1))})$ as a
direct sum of indecomposable projective $k[H]$-modules, every
summand occurs with a multiplicity divisible by $f_P$. This proves
the assertion.
\end{proof}

In the proof of Theorem \ref{my-theorem-2}(a), we have used a
preliminary version of the equivariant Riemann-Roch formula 
to show the divisibility of
$\Ind_{I_P}^{G_P} (\Cov((\frcm_P/\frcm_P^2)^{\otimes (- d)}))$ by
$f_P$, i.e.\ we have used a global argument to prove a local
statement. This tells
 us very litte about the structure of the summands $W_{P,d}$, which
 leads to the question whether one could find a ``local'' proof for
 the divisibility. In two different situations, the following
 proposition provides such a proof, yielding a concrete description
 of $W_{P,d}$.

\begin{prop} \intlabel{structure-of-WP,d}
Assume that $\pi$ is tamely ramified, let $P \in |X|$ and $d \in
\{1, \ldots, e_P^t -1\}$.
\begin{enumerate}
\item[(a)] If $\Gal(k(P)/k(\pi(P)))$ is abelian,
  then we have $W_{P,d}\cong (\frcm_{P}/\frcm_{P}^2)^{\otimes (-d)}$ as $k[G_P]$-modules.
\item[(b)] If $I_P$ is central in $G_P$, then $W_{P,d}$ is of the form
$W_{P,d}= \Ind_{I_P}^G(\chi_d)$ for some $k[I_P]$-module $\chi_d$.
If moreover $G_P \cong I_P \times G_P/I_P$, then $W_{P,d}\cong
(\frcm_{P}/\frcm_{P}^2)^{\otimes (-d)}$ as $k[G_P]$-modules.
\end{enumerate}
\end{prop}
Note that since every Galois extension of a finite field is cyclic,
the first part of this proposition gives a ``local'' proof of
Theorem \ref{my-theorem-2}(a) for the important case where $\pi$ is
tamely ramified and the underlying field $k$ is finite.

Proposition \ref{structure-of-WP,d} can be deduced from the
following purely algebraic result. Note that, in this result, we
don't use the notations introduced earlier in this paper; when
Proposition~\ref{structure-alg} is being applied to prove 
Proposition~\ref{structure-of-WP,d}, the
fields $k$ and $l$ become the fields $k(\pi(P))$ and $k(P)$,
respectively, the group $G$ becomes $G_P$ and $V$ becomes
$(\frcm_P/\frcm_P^2)^{\otimes (-d)}$ which is viewed only as a
representation of $I_P$ (and not of $G_P$) in Theorem~4.6(a).

\begin{prop} \intlabel{structure-alg}
Let $l/k$ be a finite Galois extension of fields. Let $G$ be a
finite group, and let $I$ be a cyclic normal subgroup of $G$, such
that $G/I \cong \Gal(l/k)$, i.e. we have a short exact sequence
\[ 1 \rightarrow I \rightarrow G \rightarrow \Gal(l/k) \rightarrow 1. \]
Let $V$ be a one-dimensional vector space over $l$ such that $G$
acts semilinearly on $V$, that is, for any $g \in G, \lambda \in l,
v,w \in V$, we have  $g.(\lambda v + w)= \bar g(\lambda) (g.v) +
g.w$, where $\bar g$ denotes the image of $g$ in $\Gal(l/k)$.
\begin{enumerate}
\item[(a)] If $\Gal(l/k)$ is abelian, then
we have $\Ind_I^G \Res_I^G (V) \cong \bigoplus^{(G:I)} V$ as $k[G]$-modules.
\item[(b)] If $I$ is central in $G$, then there is a (non-trivial) one-dimensional
  $k$-representation $\chi$ of $I$ such that
$\Res_I^G(V) \cong \bigoplus^{(G:I)} \chi$ as $k[I]$-modules.

If moreover $G = I \times \Gal(l/k)$, then we have $\Ind_I^G \chi \cong V$
and $\Ind_I^G \Res_I^G(V) \cong \bigoplus^{(G:I)} V$ as $k[G]$-modules.
\end{enumerate}
\end{prop}
\begin{proof}
\begin{enumerate}
\item[(a)] We have (isomorphisms of $k[G]$-modules):
\begin{align*}
\Ind_I^G & \Res_I^G (V) \\
& \cong V \otimes_k \Ind_I^G (k) \quad
 & \text{by Corollary 10.20 in \cite{CR}}\\
 & \cong V \otimes_k k[G/I] & \text{(cf. \S 10A in \cite{CR})}\\
 & \cong V \otimes_k k[\Gal(l/k)]
 & \text{as}\; \Gal(l/k) \cong G/I \\
 & \cong V \otimes_k l \\
 & \cong \bigoplus_{\sigma \in \Gal(l/k)} V.
\end{align*}
The last two isomorphisms can be derived as follows. By the normal
basis theorem, there is an element $x_0 \in l$ such that $\{g(x_0) |
g \in \Gal(l/k)\}$ is a basis of $l$ over $k$. The resulting
isomorphism
\begin{align*} k[\Gal(l/k)] & \rightarrow l \quad  \text{given by} \\
  [g] & \mapsto g(x_0) \quad \text{for every}\; g
  \in \Gal(l/k).
\end{align*} is obviously $k[G]$-linear.
This is the second last isomorphism. For the last one, we define
\begin{align*} \varphi: l \otimes_k V
& \rightarrow \bigoplus_{\sigma \in \Gal(l/k)} V
  \quad \text{by} \\
a \otimes v & \mapsto (\sigma(a) \cdot v)_{\sigma \in \Gal(l/k)}
  \quad \text{for every} \; a \in l, v \in V. \end{align*}
$\varphi$ is an isomorphism of vector spaces over $k$, by the
Galois Descent Lemma.
If $\Gal(l/k)$ is commutative, then $\varphi$ is also compatible with the
  $G$-action on both sides: Let $a \in l$, $v \in
  V$, $g \in G$, then we have
\begin{multline*} \varphi(g.(a \otimes v)) =\varphi(\bar g(a) \otimes g.v)
=((\sigma \bar g)(a) \cdot g.v)_{\sigma \in \Gal(l/k)}
= ((\bar g \sigma)(a)
  \cdot g.v)_{\sigma \in \Gal(l/k)}\\
 =  g.((\sigma(a) \cdot v)_{\sigma
  \in \Gal(l/k)})=g.\varphi(a \otimes v).\end{multline*}
\item[(b)] Since $I$ is cyclic, it acts by
multiplication with $e$-th roots of unity, where $e$ divides $|I|$.
If $I$ is central in $G$, then it follows
  that the $e$-th roots of unity are contained in $k$. For if $h$ is
  a generator of $I$ and $h.v = \zeta_e \cdot v$ for all $v \in V$,
  $\zeta_e$ 
an   $e$-th root of unity, then we have for all $g \in G$ and
  all $v \in V$:
\[ \bar g(\zeta_e)(g.v) = g.(\zeta_e v)=(gh).v=(hg).v=\zeta_e
(g.v). \]
Hence for everg $\bar g \in \Gal(l/k)$, we have
$\bar g (\zeta_e)= \zeta_e$, which means that $\zeta_e$ lies in $k$.
Let now $\{x_1, \ldots, x_f\}$ be a $k$-basis of $V$, where $f=(G:I)$.
Then we have
$V=k x_0 \oplus \ldots \oplus k x_f$ not only as vector spaces over $k$,
but also as $k[I]$-modules, since
\[I x_i = \{\zeta_e^j x_i | j=
0, \ldots , e-1\} \subseteq k x_i\] for every basis vector $x_i$.
Furthermore, the summands $k x_i$ are isomorphic as
$k[I]$-modules because $I$ acts on each of them by multiplication
with the same roots of
unity in $k$. Setting for example $k x_1 =: \chi$, we can write
\[ \Res_I^G (V) \cong \bigoplus^f \chi \]
as requested.

Assume now that $G = I \times \Gal(l/k)$. 
Then by the Galois Descent Lemma, we have
\[V \cong l \otimes_k V^{\Gal(l/k)}\]
as $k[G]$-modules, where $I$ acts trivially on $l$ and $\Gal(l/k)$
acts trivially on $V^{\Gal(l/k)}$.
This is isomorphic to $l \otimes_k \chi$, where  $\chi$ is regarded as a
$k[G]$-module via the projection $G= I \times \Gal(l/k) \rightarrow
I$. By the normal basis theorem, we have
\[l \otimes_k \chi \cong \Ind_I^G (k) \otimes \chi = \Ind_I^G
(\chi),\]
so $V \cong \Ind_I^G(\chi)$ as requested. Together with what we have
shown before,
this implies the last identity of the proposition:
\[\Ind_I^G \Res_I^G (V) = \Ind_I^G (\bigoplus^f \chi) = \bigoplus^f
V.\]
\end{enumerate}
\end{proof}

\section{Some variants of the main theorem}
Throughout the previous section, we have concentrated on the case
where $\pi: X \rightarrow Y$ is weakly ramified and where the locally free
$G$-sheaf we are considering comes from an equivariant divisor. 
If $\pi$ is tamely ramified, we have the
following variant of Theorem~\ref{my-theorem-2} for locally free
$G$-sheaves that need not come from a divisor. It generalizes
Corollary 1.4(b) in \cite{BK1}.
\begin{theorem} \label{theorem-for-sheaves}
Let $\pi: X \rightarrow Y$ be tamely ramified. Let $\scrE$ be a
locally free $G$-sheaf of rank $r$ on $X$. 
For every closed point $P \in |X|$ and for $i=1, \ldots, r$, let the integers 
$l_{P,i} \in \{0, \ldots , e_P -1\}$ be defined by the following
isomorphism of $k(P)[I_P]$-modules:
\[\scrE(P) \cong \bigoplus_{i=1}^r \bigm(\frcm_P/\frcm_P^2 \bigm)^{\otimes l_{P,i}}.\]
For every $R \in |Y|$, let $\tilde R \in |X|$ and $W_{\tilde R,d}$
be defined as in Theorem \ref{my-theorem-2}. 
Furthermore, let $N_{G,X}$ be the ramification
module from Theorem \ref{NG,X}. 
Then we have in $K_0(k[G])$:
\[ \chi(G,X,\scrE) \equiv -r[N_{G,X}] + \sum_{R \in Y} \sum_{i=1}^r
\sum_{d=1}^{l_{\tilde R, i}} [\Ind_{G_{\tilde R}}^G(W_{\tilde R,d})]
\mod \ZZ[G]. \]
\end{theorem}

Moreover, one can show an equivariant Riemann-Roch formula
for \emph{arbitrarily ramified covers} $\pi: X \rightarrow Y$. Recall that 
in Theorem \ref{Euler-char-projective}, we have shown that in virtually
all cases where the Euler
characteristic lies in the image of the Cartan homomorphism, the cover
$\pi$ is weakly ramified. So in the general case, 
one cannot possibly find a formula in the
Grothendieck group $K_0(k[G])$ of projective $k[G]$-modules. 
However, in the Grothendieck
group $K_0(G,k)$ of all $k[G]$-modules, we have the following result, which
generalizes Theorem 3.1 in \cite{BK2}.

\begin{theorem} \label{my-theorem}
Let $\scrE$ be a locally free $G$-sheaf. Then we have in $K_0(G,k)$:
\[ n \, \chi(G,X,\scrE) = C_{G,X,\scrE}\; [k[G]] - \sum_{P \in |X|} e_P^w
\sum_{d=0}^{e_P^t -1} d\, [\Ind_{I_P}^G \bigm(\scrE(P) \otimes_{k(P)}
(\frcm_P/\frcm_P^2)^{\otimes d}\bigm) ], \]
where
\[ C_{G, X, \scrE}= r (1 - g_X) + \deg \scrE + \frac{r}{2} \sum_{P \in
|X|} [k(P):k] (e_P^t -1).\]
\end{theorem}

We omit the proofs of Theorem \ref{theorem-for-sheaves} and Theorem
\ref{my-theorem} due to their similarity with the proof of Theorem 
\ref{my-theorem-2}.

\end{document}